\documentclass[12pt]{article}

\usepackage{amsmath,amssymb,mathtools,amsthm}
\mathtoolsset{showonlyrefs}
\usepackage{dsfont}
\usepackage{enumerate}
\usepackage{calc}
\usepackage{authblk}
\usepackage[utf8]{inputenc}
\usepackage[
    inner=2.3cm,
    outer=2.3cm,
    top=2.7cm,
    marginparwidth=2.5cm,
    marginparsep=0.1cm
    ]
    {geometry}
\usepackage[dvipsnames]{xcolor}
\usepackage{overpic}
\usepackage{color}
\usepackage[pdffitwindow=false,
            plainpages=false,
            pdfpagelabels=true,
            pdfpagemode=UseOutlines,
            pdfpagelayout=SinglePage,
            bookmarks=false,
            colorlinks=true,
            hyperfootnotes=false,
            linkcolor=blue,
            urlcolor=blue!30!black,
            citecolor=green!50!black]{hyperref}
\usepackage{appendix}
\newcommand{\R}{\mathds{R}}
\newcommand{\N}{\mathds{N}}
\newcommand{\C}{\mathds{C}}
\newcommand{\Z}{\mathds{Z}}
\newcommand{\E}{\mathds{E}}
\newcommand{\norm}[1]{\left\lVert#1\right\rVert}
\newcommand{\abs}[1]{\left|#1\right|}

\renewcommand{\P}{\mathds{P}}
\newcommand{\tor}{\mathds{T}}
\newcommand{\Id}{\textup{Id}}

\newcommand{\dd}{\mathrm d}

\newcommand{\revision}[1]{#1}

\newenvironment{nalign}{
    \begin{equation}
    \begin{aligned}
}{
    \end{aligned}
    \end{equation}
    \ignorespacesafterend
}

\newtheorem{theorem}{Theorem}[section]
\newtheorem{proposition}[theorem]{Proposition}
\newtheorem{lemma}[theorem]{Lemma}
\newtheorem{corollary}[theorem]{Corollary}

\newtheorem{conjecture}{Conjecture}

\theoremstyle{definition}
\newtheorem{remark}[theorem]{Remark}

\begin{document}
\title{Mixing at the Batchelor Scale for White-In-Time Flows}
\author[1]{Robin Chemnitz\thanks{Corresponding author. Email: \texttt{r.chemnitz@fu-berlin.de}}}
\author[2]{Dennis Chemnitz\thanks{Email: \texttt{dennis.chemnitz@math.ethz.ch}}}

\affil[1]{Freie Universität Berlin}
\affil[2]{ETH Zürich}

\maketitle
\begin{abstract}
    We consider the mixing properties of solutions to the advection-diffusion equation of a white-in-time velocity field on the 2-dimensional torus with four forced modes. As the diffusivity parameter goes to zero, we show that the almost-sure exponential dissipation rate stays bounded from below. Together with the corresponding upper bound established by Gess and Yaroslavtsev, this constitutes an example of a velocity field for which the Batchelor scale conjecture can be verified. In addition, we characterize the exponential mixing rate without diffusion of this system. Our results are not restricted to two dimensions, and we construct a three-dimensional white-in-time velocity field with the same properties.
\end{abstract}
\section{Introduction}

In this paper, we discuss the mixing properties of passive particles on the two-dimensional torus advected by a white-in-time shear flow. 
In particular, we study the long-term behavior of the SPDE
\begin{equation}\label{spde:main2D}
    \dd f_t = \kappa \Delta f_t \, \dd t - \sum_{j = 1}^4 \sigma_j \cdot \nabla f_t \circ \dd W_t^j,
\end{equation}
where $f : [0,\infty) \times \mathbb T^2 \to \mathbb R$ is the concentration of particles initialized at $f_0 \in L^2(\mathbb T^2)$, $0 \leq \kappa \ll 1$ is a real parameter describing the strength of molecular diffusion, $W^1, \dots, W^4$ are independent Brownian motions, and  $\sigma_1, \dots, \sigma_4$ are the divergence-free vector fields given by
\begin{equation} \label{eq:4modes}
    \begin{aligned}
        \sigma_1(x,y) &:= (\sin(2\pi y),0)^\top,&\ \sigma_2(x,y) &:= (\cos(2\pi y),0)^\top\\
    \sigma_3(x,y) &:= (0,\sin(2\pi x))^\top,& \sigma_4(x,y) &:= (0,\cos(2\pi x))^\top.
    \end{aligned}
\end{equation}
Our main result, Theorem \ref{thm:lower_bound_low_modes} below, states that the $L^2$-mass of $f_t$ in the four lowest Fourier modes can decay at most exponentially with a rate that is bounded away from $-\infty$ uniformly in $\kappa$. Together with the enhanced dissipation \revision{results from \cite{gess2025stabilization, luo2024elementary}}, this implies that the dissipation rate $\lambda^\kappa$ is of constant order for small $\kappa$ and the typical filamentation length $\ell(f_t) \coloneqq \norm{ f_t}_{L^2} / \norm{\nabla f_t}_{L^2}$ obeys the Batchelor scale $\ell \sim \sqrt{\kappa}$ \cite{miles2018diffusion}; see Corollary 
\ref{cor:BatchScale} and Figure \ref{fig:examples} below. To our knowledge this is the first model of passive scalar \revision{advection} for which this can be shown. In the case of $\kappa=0$, Theorem \ref{thm:lower_bound_low_modes} allows us to determine the dependence of the $H^{s} \to H^{-s}$-mixing rate $\gamma_s$ on $s$ up to a constant; see Corollary \ref{cor:mixingRate} below. Our method works in higher dimensions, and we construct a three-dimensional model analogous to \eqref{spde:main2D} for which the same results hold; see Section \ref{sec:3d}.

\newcommand{\tickcolor}{\color{red}}
\begin{figure}
    \centering

    \includegraphics[width=0.8\linewidth]{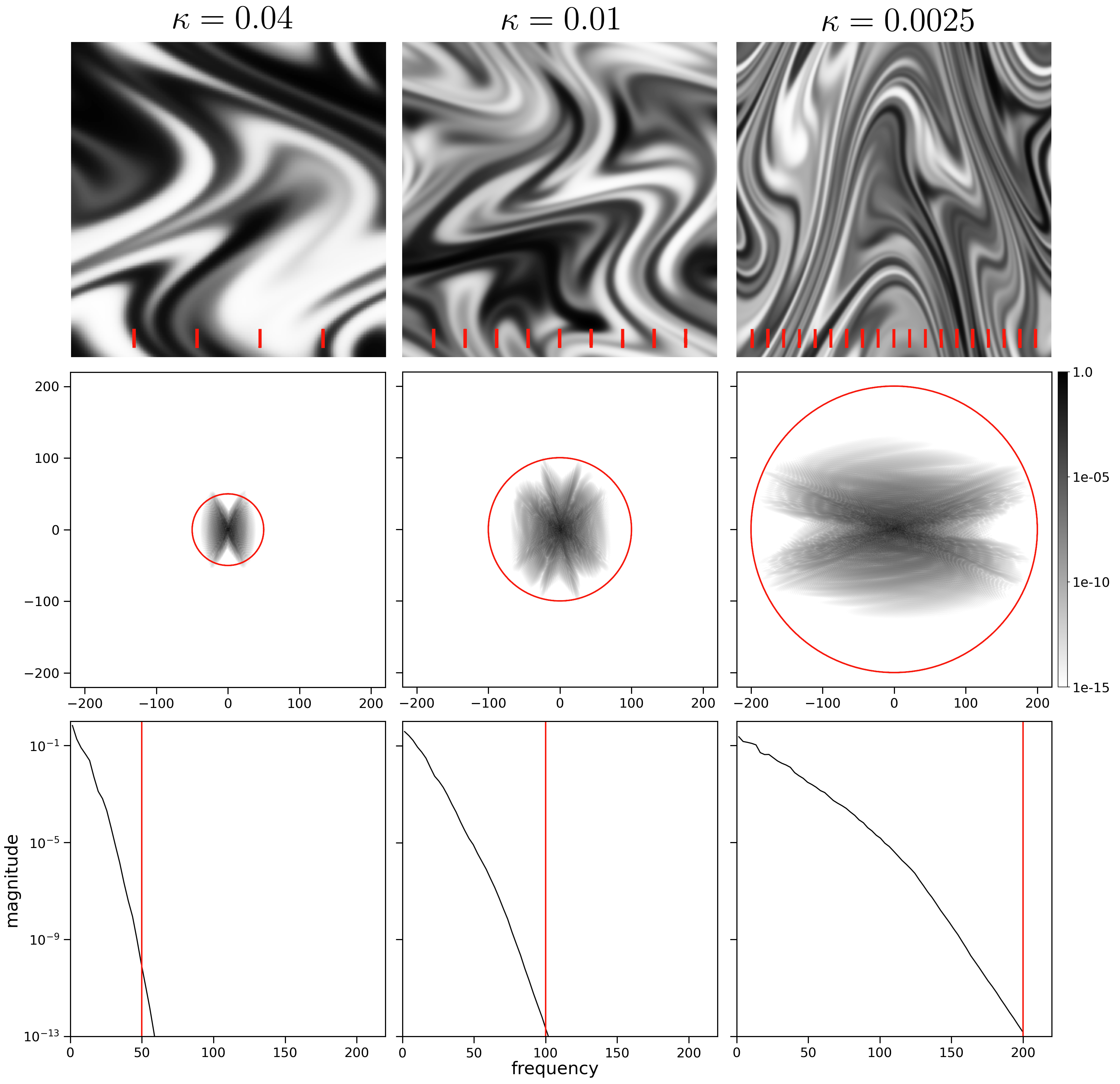}
    \caption{Each column represents a typical realization of the normalized solution $\pi_t \coloneqq f_t / \norm{f_t}_{L^2}$ to the 4-modes model \eqref{spde:main2D} for different values of $\kappa$; from left to right $\kappa=0.04, \, 0.01,\, 0.0025$. The first row depicts the function $\pi_t$ on the torus, and the predicted filamentation length $\sqrt{\kappa}$ is marked in red. The second row shows the Fourier spectrum of $\pi_t$ as a heatmap on $\Z^2$. A multiple of the Batchelor scale ($10\,\kappa^{-\frac{1}{2}}$) is drawn as a red circle. The third row depicts the power spectrum of $\pi_t$ with a red line at $10\, \kappa^{-\frac{1}{2}}$. \revision{The numerical simulations use a Galerkin projection of the SPDE \eqref{spde:main2D} onto the Fourier modes $[-10^3, 10^3] \times[-10^3, 10^3]$. In each time step of length $\Delta t = 10^{-2}$, the increments $\dd W_t^1, \hdots, \dd W_t^4$ are drawn, and $f_{t+\Delta t}$ is computed using a Lie--Trotter splitting of the advection and diffusion operators. }
    }
    \label{fig:examples}
\end{figure}

\subsection{Dissipation of a passive scalar in the Batchelor regime}
 Understanding the nature of turbulence occurring in three-dimensional high Reynolds number fluid flows is one of the major open problems in continuum mechanics. In particular, the Kolmogorov $-5/3$ scaling, which predicts the energy distribution of a turbulent fluid flow over Fourier modes \cite{kolmogorov1991local,frisch1995turbulence}, is wide open. A related, but mathematically more tractable, problem concerns the Fourier mass distribution of the solution to an advection-diffusion equation for which the advection vector field is given \revision{by a stochastic fluid model}. \revision{Such an equation takes the form
 \begin{equation}\label{pde:advdiff}
    \partial_t f_t = \kappa \Delta f_t - u_t \cdot \nabla f_t,
\end{equation}
where $u$ is a random, time-dependent, divergence-free velocity field and $f:[0,\infty) \to L^2(\mathbb T^d)$ models the concentration of passive particles dissolved in the fluid.}
 Depending on the ratio of the Reynolds number of $u$ and the diffusivity $\kappa$, different \revision{behavior is expected.}

\revision{In this work, we consider the so-called Batchelor regime, in which there has been tremendous progress in the understanding of passive scalar mixing over the last years \cite{bedrossian2021PTRF,bedrossian2022JEMS, bedrossian2022AOP, bedrossian2022CPAM, hairer2024lower,cooperman2025fourier,  gess2025stabilization}.}
The Batchelor regime \cite{batchelor1959small} refers to a setting where the diffusivity $\kappa$ is small while the Reynolds number of the fluid is fixed and finite.
Typical examples of such random velocities $u$ are the solutions to the 2- and 3-dimensional stochastic (hyperviscous) Navier--Stokes equation with positive viscosity $\nu>0$ and adequately chosen forcing (see, e.g., \cite[Systems 1 and 2]{bedrossian2021PTRF} for details).
\revision{For velocity fields $u$ which are not given as solutions to the Navier--Stokes equation, the assumption of a finite Reynolds number can be replaced by assuming sufficient spatial smoothness of $u$, namely  $u \in C([0,\infty),C^1(\mathbb T^d, \mathbb R^d))$.}

For simplicity, we assume that $f_0 \in L^2_0$ is mean-free, where $L^2_0:= \{f \in L^2(\mathbb T^d):\int_{\mathbb T^d} f(x)\, \dd x = 0\}$. By conservation of mass, this implies that $f_t$ is mean-free for every $t\geq 0$. The $L^2$-norm of $f_t$ serves as a measure of the distance between the concentration of particles $f_t$ and the equilibrium. Consequently, the exponential dissipation rate can be measured by
\begin{equation}\label{def:LE1}
    \lambda^{\kappa}_{u, f_0} := \limsup_{t\to \infty} \frac{1}{t} \log \|f_t\|_{L^2}.
\end{equation}

In the following, we assume that $u$ is random and distributed  according to some ergodic stationary measure $\mu$ of divergence-free velocity fields. Let $\mathcal T_u^t: f_0 \mapsto f_t$ be the time-$t$ operator of the PDE \eqref{pde:advdiff}, which, by classical PDE theory, is a compact linear operator from $L_0^2(\mathbb T^d)$ to itself with $\|\mathcal T_u^t\|_{L_0^2} \leq 1$. By the Multiplicative Ergodic Theorem \cite{oseledec1968multiplicative, ruelle1982characteristic}, the quantity $\lambda^{\kappa}_{u, f_0}$ defined in \eqref{def:LE1} is
a true limit. 
Additionally, the supremum of $\lambda^\kappa_{u, f_0}$ over all initial conditions is essentially independent of $u$. There exists a $\lambda^\kappa$ such that for almost every $u$ with respect to $\mu$ we have
$$\lambda^{\kappa} = \sup_{f_0 \in L^2_0} \lambda^{\kappa}_{u, f_0} = \lim_{t \to \infty} \frac{1}{t} \log\|\mathcal T^t_{u}\|_{L^2_0} \in [-\infty, 0].$$
By a standard calculation, it can be shown that
\begin{align}\label{eq:ell}
    \frac{\dd}{\dd t}\|f_t\|_{L^2}^2 = -2\kappa \|\nabla f_t\|_{L^2}^2 \quad \Rightarrow\quad \frac{\dd}{\dd t} \log \|f_t\|_{L^2} = - \kappa \frac{\|\nabla f_t\|^2_{L^2}}{\|f_t\|^2_{L^2}} \eqqcolon -\kappa \ell(f_t)^{-2}.
\end{align}
The quantity $\ell(f_t) \coloneqq \|f_t\|_{L^2}/\|\nabla f_t\|_{L^2}$ can be interpreted as a typical ``filamentation length" of the passive scalar concentration \cite[Section 1.2]{hairer2024lower}.
By \eqref{eq:ell}, the filamentation length $\ell(f_t)$ is closely related to the dissipation rate $\lambda^{\kappa}_{u,f_0}$ through
\begin{equation}\label{eq:lambdaEll}
    \lambda^{\kappa}_{u, f_0} = \limsup_{t \to \infty} -\kappa\frac{1}{t} \int_0^t \ell(f_s)^{-2}\,\dd s.
\end{equation}
The Poincaré inequality gives upper bounds $\ell(f_t) \lesssim 1$ for the filamentation length and ${\lambda^{\kappa}_{u, f_0}\lesssim -\kappa}$ for the dissipation rate, with a constant that is independent of $\kappa, u$ and $f_0$. These bounds are optimal for $u = 0$. However, for velocity fields $u$ with chaotic Lagrangian trajectories, one expects the filamentation to become finer for a smaller diffusion coefficient $\kappa$. In fact, numerical evidence \cite{miles2018diffusion}; see also Figure \ref{fig:examples}, suggests that the filamentation length decays proportionally to the diffusive length scale $\ell(f_t)\sim \sqrt \kappa$. By the relation \eqref{eq:lambdaEll}, this corresponds to a constant order dissipation rate $\lambda^\kappa_{f_0} \sim -1$ for small $\kappa$.

\begin{conjecture}[Batchelor Scale Conjecture]\label{conj:Batchelor}
For $u$ given by a 2-dimensional stochastic Navier--Stokes equation \cite[Systems 1]{bedrossian2021PTRF} or a 3-dimensional hyperviscous Navier--Stokes equation \cite[Systems 2]{bedrossian2021PTRF}, the dissipation rate $\lambda^\kappa$ is of order 1 for small $\kappa$, that is
$$ -\infty<\liminf_{\kappa \to 0} \lambda^\kappa \leq \limsup_{\kappa \to 0} \lambda^\kappa < 0.$$
\end{conjecture}
The upper bound $\lambda^\kappa_{u, f_0} \lesssim -1$ has been established in the seminal work \cite{bedrossian2021PTRF} by Bedrossian, Blumenthal, and Punshon-Smith.
Somewhat surprisingly, the corresponding lower bound is wide open. In fact, even showing that $\lambda^\kappa > -\infty$ for any fixed $\kappa$ requires significant work. For the 2-dimensional stochastic Navier--Stokes equation, Hairer, \revision{Punshon-Smith}, Rosati, and Yi \cite{hairer2024lower} showed that $\lambda^\kappa \gtrsim \kappa^{-q}$ for any $q>3$. This constitutes the first proof of a finite lower bound for the dissipation rate of any non-autonomous velocity field. (For general autonomous vector fields, the situation is easier and the lower bound $\lambda^\kappa > - \infty$ can be established using PDE techniques, for example \cite[Theorem 16.4]{agmon2010lectures}.)

\revision{While Rowan~\cite{rowan2025superexponential} constructed examples of deterministic vector fields with superexponential dissipation for certain initial conditions, there is to date no example of a smooth velocity field with $\lambda_\kappa = -\infty$. For discrete-time dynamical systems with pulsed diffusion, that is, alternating a transport operator with a diffusion operator of strength $\kappa$, there are multiple known examples for which $\lambda_\kappa=-\infty$, see \cite{feng2019dissipation, gonzalez2021stability}. This is known as \emph{collapse of the Lyapunov spectrum}.
}

\subsection{Almost-sure exponential mixing in the Batchelor regime}\label{sec:exponential_mixing}
\revision{Another quantity that is closely related to the diffusion rate is the so-called \emph{mixing rate.}}
For $\kappa = 0$, equation \eqref{pde:advdiff} turns into the transport equation
\begin{equation}\label{pde:Transport}
    \partial_t f_t = - u_t\cdot \nabla f_t.
\end{equation}
Since the vector field $u_t$ is divergence-free and sufficiently regular, the $L^2$-norms of (weak) solutions are conserved, and we have $\lambda_{u,f_0}^0 = 0$ for every $f_0$. In this setting, the mixing of passive scalars can be measured in terms of the evolution of a negative regularity Sobolev norm $\|f_t\|_{H^{-s}}$ with $s>0$; see \eqref{eq:Sobolev_norm}.

The transport equation \eqref{pde:Transport} is said to be almost-surely \emph{exponentially mixing} \cite{bedrossian2022AOP} if for some $s>0$, there exist $\gamma>0$ and a constant depending on the vector field $C_{u}>0$, such that for almost every $u$
\begin{equation}\label{ineq:mixing}
   \|f_t\|_{H^{-s}} \leq C_{u} e^{-\gamma t} \|f_0\|_{H^s},\quad\forall\, f_0 \in H^s_0 \coloneqq H^s\cap L^2_0. 
\end{equation}
Establishing exponential mixing for random dynamical systems has been a topic of active research \cite{dolgopyat2004sample, blumenthal2023exponential, dewitt2024expanding, navarro2025exponential}. 
In the context of fluid dynamics, Bedrossian, Blumenthal, and Punshon-Smith showed that the solutions of the 2-dimensional stochastic Navier--Stokes equation and the 3-dimensional hyperviscous stochastic Navier--Stokes equation are almost surely exponentially mixing \cite{bedrossian2022AOP}; see also \cite{cooperman2025exponential} for the case of degenerate noise.

Given a velocity field, let $\gamma_s$ be the supremum of all values of $\gamma$ for which \eqref{ineq:mixing} holds, where the constant $C_u$ may depend on $s$. If there is no such $\gamma$, we set $\gamma_s = 0$. We call $\gamma_s$ the \emph{mixing rate} for the regularity parameter $s$. Thus, a system is exponentially mixing if there is an $s>0$ such that $\gamma_s>0$. Without any assumptions on the equation \eqref{pde:Transport}, it can be shown  that if $\gamma_{s_0}>0$ for some $s_0>0$, then $\gamma_s>0$ for all $s>0$ (see Proposition \ref{prop:GammaSLowerBound}, also \cite[Section 7]{bedrossian2022AOP}). In particular, if some mixing rate $\gamma_{s_0}$ is known, the mixing rate $\gamma_s$ for any $s$ can be lower bounded by
$$\gamma_s \geq \frac{s}{2s_0-s}\gamma_{s_0}, \quad \text{for all } 0<s<s_0\quad \quad\text{and}\quad \quad\gamma_s \geq \gamma_{s_0}\quad\text{for all } s>s_0.$$
Consequently, if \eqref{pde:Transport} is exponentially mixing, then $\gamma_s \gtrsim 1 \wedge s$.

For velocity fields with a Lipschitz constant that is either uniformly bounded in time, or satisfies an appropriate moment bound, there is an upper bound on the mixing rate $\gamma_s$ that is linear in $s$. For $u$ given by the solution to the stochastic Navier--Stokes equation, this upper bound can be shown using results from \cite[Lemma 7.3 and equation (7.2)]{bedrossian2021PTRF} and an interpolation argument; see Lemma \ref{lemm:interpolation}. 

Hence, for exponentially mixing velocity fields with an appropriate moment bound, there is a constant $C>0$ such that for all $s>0$
\begin{equation*}
    C^{-1} (1 \land s) \leq \gamma_s \leq Cs.
\end{equation*}
For large $s$, it is not clear whether $\gamma_s$ continues to grow linearly, whether it grows sublinearly, or whether it remains bounded (cf.~Corollary \ref{cor:mixingRate}) Our main results, which we present below, suggest a deeper connection between this question and Conjecture \ref{conj:Batchelor}.

\revision{While for pure transport equations ($\kappa=0$) the mixing rate $\gamma_s$ depends on the parameter $s$, the exponential dissipation rate of an advection-diffusion equation ($\kappa>0$) is less sensitive to the norm chosen.
Consider the dissipation rate
\begin{equation}
    \lambda_{u, f_0}^{\kappa, s} \coloneqq \limsup_{t\to \infty} \frac{1}{t} \log \norm{f_t}_{H^s},
\end{equation}
measured in a Sobolev norm $\norm{\cdot}_{H^s}$. A result by Blumenthal and Punshon-Smith \cite[Theorem 4.2]{blumenthal2023norm} states that for sufficiently smooth (in space) velocity fields, $\lambda_{u, f_0}^{\kappa, s}$
is the same for all $s\in [-K, K]$, where $K>0$ depends on the spatial smoothness of the velocity field.}

\section{White-in-time velocity fields and main results}\label{sec:white}

An aspect that makes the analysis of $f_t$ evolving under \eqref{pde:advdiff} challenging is the fact that the velocity field $u_t$ has a memory. This causes the process $f_t$ to be non-Markovian. To avoid these difficulties, one can consider a different model for the evolution of passive scalars, given by SPDEs of the form
\begin{equation}\label{spde:advDiff}
    \dd f_t = \kappa \Delta f_t \, \dd t - \sum_{j = 1}^m \sigma_j \cdot \nabla f_t \circ \dd W_t^j.
\end{equation}
Here $m \in \mathbb N$, $\alpha \in (0,1]$, $\sigma_1,\dots, \sigma_m \in C^{2, \alpha}(\mathbb T^d, \mathbb R^d)$ are autonomous, divergence-free vector fields on the torus, and $W^1, \dots, W^m$ are independent Brownian motions. 
\revision{The notation $\circ \,\dd W_t^j$ indicates that the corresponding stochastic integrals should be interpreted in the Stratonovich sense. Unlike the It\^o convention, the Stratonovich advection terms with sufficiently smooth coefficients preserve the $L^2$-norm of solutions. Furthermore, by the Wong--Zakai theorem \cite{BRZEZNIAK1995329}, the SPDE \eqref{spde:advDiff} can be seen as an approximation for systems of the form \eqref{pde:advdiff} where the velocity field has a short correlation time \cite{kramer2003two}. In this analogy, the velocity field $u$ is formally given by 
$$u_t(x) = \sum_{j = 1}^m \sigma_j(x) \dot W^j_t.$$
}


Again, we consider mean-free initial conditions $f_0 \in L^2_0.$ Analogous to the introduction, let $\mathcal T_\omega^t: f_0 \mapsto f_t$ be the random time-$t$ operator of the SPDE \eqref{spde:advDiff}. As before, \revision{$\mathcal T_\omega^t: L_0^2(\mathbb T^d) \to L_0^2(\mathbb T^d)$ is almost surely compact for every $\kappa, t>0$ (see Appendix \ref{app:compact}) and thus} the Multiplicative Ergodic Theorem \cite{ruelle1982characteristic} implies that for each $\kappa>0$ and $f_0 \in L^2_0$ the exponential dissipation rate $\lambda_{\omega, f_0}^\kappa$ is well-defined for almost every $\omega$,
\begin{equation}\label{eq:dissipation_rate}
    \lambda_{\omega, f_0}^\kappa \coloneqq \lim_{t \to \infty} \frac{1}{t} \log \|\mathcal T_\omega^t f_0\|_{L^2}.
\end{equation}
Additionally, the supremum over all initial conditions essentially does not depend on $\omega$. There exists a $\lambda^\kappa$ such that for almost every $\omega$ we have 
\begin{equation}
    \lambda^\kappa \coloneqq \sup_{f_0 \in L^2_0} \lambda^{\kappa}_{\omega, f_0} = \lim_{t \to \infty} \frac{1}{t} \log\|\mathcal T^t_{\omega}\|_{L^2_0} \in [-\infty, 0].
\end{equation}
Gess and \revision{Yaroslavtsev} \cite{gess2025stabilization}\revision{, as well as Luo, Tang and Zhao \cite{luo2024elementary}}, showed that under mild conditions on the functions $\sigma_i$, the dissipation rate $\lambda^\kappa$ has a negative upper bound for small $\kappa$, that is, 
\begin{equation}\label{eq:GY_upper_bound}
    \limsup_{\kappa \to 0} \lambda^\kappa < 0.
\end{equation}
In particular, this upper bound holds for the 4-modes model \cite[Theorem 6.2]{gess2025stabilization}, where $m =4$ and $\sigma_1,\dots, \sigma_4$ given by \eqref{eq:4modes}.

\subsection{Main Results}

 The main contribution (Corollary \ref{cor:BatchScale} below) of this paper is the corresponding lower bound 
\begin{equation}\label{ineq:LowerBound}
    \liminf_{\kappa \to 0} \lambda^\kappa > -\infty,
\end{equation}
in the 4-modes model \eqref{spde:main2D}. To our knowledge, this constitutes the first example of a model for passive scalar mixing in which the Batchelor scale conjecture can be verified. Our main result Theorem \ref{thm:lower_bound_low_modes} is formulated in terms of Fourier coefficients, which we briefly introduce to clarify notation.

For scalar functions on the $d$-dimensional torus, Fourier modes are indexed by $\mathbf{k} \in \Z^d$ and given by $F_{\mathbf{k}} (\mathbf{x}) = e^{2\pi i\,\mathbf{k} \cdot \mathbf{x}}$. The Fourier coefficients of $f_t$ are denoted by
\begin{equation}
    \widehat{f}_{\mathbf{k}}(t) = \langle F_\mathbf{k}, f_t\rangle =  \int_{\tor^d} e^{-2\pi i\,\mathbf{k} \cdot \mathbf{x}} \,f_t (\mathbf{x}) \, \dd \mathbf{x}.
\end{equation}
Since $f_t$ is assumed to be mean-free and real-valued, we have $\widehat{f}_\mathbf{0}(t) = 0$ and $\widehat{f}_{-\mathbf{k}}(t) = \overline{\widehat{f}_\mathbf{k}(t)}$. For a function $f \in L_0^2(\mathbb T^d)$ the Sobolev norm with regularity parameter $s\in \R$ can be defined by
\begin{equation} \label{eq:Sobolev_norm}
    \|f\|^2_{H^{s}} \coloneqq \sum_{\mathbf{k}  \in \mathbb Z^d\setminus \{0\}} |\widehat f_{\mathbf{k}}|^2\,\revision{|\mathbf{k}|^{2s}}.
\end{equation}
 
The heart of our argument is to control the Fourier mass in the four innermost modes of $f_t$, that is, we consider the projection

\begin{equation}
    \Pi_{\leq 1} f_t \coloneqq \widehat{f}_{0,1}(t) F_{0,1} + \widehat{f}_{0,-1}(t) F_{0,-1} + \widehat{f}_{1,0}(t) F_{1,0} + \widehat{f}_{-1,0}(t)  F_{-1,0}.
\end{equation}
The following theorem provides a lower bound for the $L^2$-norm of this projection.

\begin{theorem}\label{thm:lower_bound_low_modes}
    For $\kappa\geq 0$ and $f_0\in L_0^2(\tor^2) \setminus \{0\}$, let $f_t$ be the solution to \eqref{spde:main2D}. The decay rate of the $L^2$-norm of the four innermost modes is almost surely lower bounded by
    \begin{equation}\label{ineq:MainLowerBound}
        \limsup_{t\to \infty} \frac{1}{t} \log \norm{\Pi_{\leq 1}f_t}_{L^2} \geq -(2\pi)^2 \left( \frac{1}{2} + \kappa \right ).
    \end{equation}
\end{theorem}
The proof of Theorem \ref{thm:lower_bound_low_modes} is given in Section \ref{sec:proof}.
As an immediate consequence, using the estimate $\|f_t\|_{L^2} \geq \|\Pi_{\leq 1}f_t\|_{L^2}$, we obtain the claimed lower bound for $\lambda^\kappa_{\omega, f_0}$.
\begin{corollary} \label{cor:BatchScale}
    For all $\kappa > 0$ and $f_0 \in L_0^2(\tor^2) \setminus \{0\}$, the exponential dissipation rate $\lambda_{\omega, f_0}^\kappa$, as given by \eqref{eq:dissipation_rate}, is almost surely lower bounded by 
    \begin{equation}
        \lambda_{\omega, f_0}^\kappa  \geq -(2\pi)^2 \left( \frac{1}{2} + \kappa \right ).
    \end{equation}
\end{corollary}
Together with the upper bound \eqref{eq:GY_upper_bound}, which was established in \cite[Theorem 6.2]{gess2025stabilization}, Corollary \ref{cor:BatchScale} verifies the Batchelor scale conjecture for the 4-modes model \eqref{spde:main2D}. 

Another interesting aspect of Theorem \ref{thm:lower_bound_low_modes} is that the lower bound holds even for the pure-advection case $\kappa = 0$. In this case, equation \eqref{spde:main2D} turns into the stochastic transport equation
\begin{equation}\label{spde:Transport}
    \dd f_t = - \sum_{j = 1}^4 \sigma_j \cdot \nabla f_t \circ \dd W_t^j.
\end{equation}
For $s>0$, let $\gamma_s$ be the mixing rate of the transport equation, as defined in Section \ref{sec:exponential_mixing}. That is, $\gamma_s$ is the supremum over all $\gamma$ for which there is a random constant $C_\omega>0$ that is almost surely finite such that
\begin{equation}\label{ineq:random_mixing}
   \|f_t\|_{H^{-s}} \leq C_{\omega} e^{-\gamma t} \|f_0\|_{H^s},\quad\forall\, f_0 \in H^s_0 \coloneqq H^s\cap L^2_0. 
\end{equation}
In \cite{gess2025stabilization} it is shown that the transport equation \eqref{spde:Transport} is exponentially mixing for any $s>0$. The following corollary shows that Theorem \ref{thm:lower_bound_low_modes} allows us to determine the dependence of $\gamma_s$ on $s$ up to a constant. 
\begin{corollary}\label{cor:mixingRate}
    We have $\gamma_s \sim 1\wedge s$, that is, there exists a constant $C>0$ such that
    \begin{equation}\label{eq:gamma_s_scaling}
        C^{-1} (1\wedge s) \leq \gamma_s \leq C (1\wedge s), \quad \forall\, s>0.
    \end{equation}
\end{corollary}
\begin{proof}
    The proof consists of four inequalities. We show that there are constants $C_1,c_2, C_3, c_4$ such that
    \begin{enumerate}[(i)]
        \item $\gamma_s \leq C_1s$ for $s\in (0, 1)$; 
        \item $\gamma_s \geq c_2s$ for $s\in (0, 1)$; 
        \item $\gamma_s \leq C_3$ for $s\geq 1$; 
        \item $\gamma_s \geq c_4$ for $s\geq 1$. 
    \end{enumerate}
    Once these constants are found, the estimate \eqref{eq:gamma_s_scaling} holds with $C\coloneqq \max\{C_1, c_2^{-1}, C_3, c_4^{-1}\}$.
    
    The first inequality is shown in Proposition \ref{prop:GammaSUpperBound}. The second inequality follows from \cite[Theorem 6.2]{gess2025stabilization}. To prove the third inequality, observe that for any $s>0$ the $H^{-s}$-norm of a function $f$ is lower bounded by $\norm{\Pi_{\leq 1} f}_{L^2}$; see \eqref{eq:Sobolev_norm}. Our main result Theorem \ref{thm:lower_bound_low_modes} shows
    $$ \limsup_{t\to \infty} \frac{1}{t} \log\|f_t\|_{H^{-s}} \geq \limsup_{t\to \infty} \frac{1}{t} \log \norm{\Pi_{\leq 1}f_t}_{L^2} \geq -2\pi^2,$$
    which implies the upper bound $\gamma_s \leq 2\pi^2$. The fourth inequality is immediate, since $\gamma_s \geq \gamma_1$ for all $s\geq 1$. This completes the proof.
\end{proof}
\begin{remark}
    \revision{The proof of Theorem \ref{thm:lower_bound_low_modes} relies on the specific structure of the 4-modes model \eqref{spde:main2D}. While the strategy of the proof is applicable to more general models of the form \eqref{spde:advDiff}, the core of our arguments consists of checking the algebraic condition \eqref{eq:algebraic_condition}. We did not find any other two-dimensional model of the form \eqref{spde:advDiff} for which we are able to verify this condition. The crucial property of the 4-modes model is that in the Fourier representation of the SPDE, cf.~Proposition \ref{prop:Fourier_representation} and Figure \ref{fig:modes}, only neighbouring Fourier modes directly interact with one another. In the following section, we consider an analogous model in three dimensions whose velocity fields consist of the lowest 12 Fourier modes. For this model, the algebraic condition \eqref{eq:algebraic_condition} is verifiable, and we obtain results analogous to Theorem \ref{thm:lower_bound_low_modes}.}
\end{remark}

\subsection{Three-dimensional model}\label{sec:3d}

In the \revision{physics} community, there is a natural interest in three-dimensional fluid models. We show that our approach works in three dimensions and we establish the same results as for the two-dimensional model \eqref{spde:main2D}. The analogous three-dimensional model is given by
\begin{equation}\label{spde:3D}
    \dd f_t = \kappa \Delta f_t \, \dd t - \sum_{j = 1}^{12} \sigma_j \cdot \nabla f_t \circ \dd W_t^j,
\end{equation}
where $f_t : [0, \infty) \times \tor^3 \to \R$, and $\sigma_1, \dots, \sigma_{12}$ are the shear flows given by
\begin{equation} \label{eq:12modes}
    \begin{aligned}
        \sigma_1(x,y,z) &:= (\sin(2\pi y),0, 0)^\top,&\ \sigma_2(x,y,z) &:= (\cos(2\pi y),0, 0)^\top\\
        \sigma_3(x,y, z) &:= (0, \sin(2\pi x), 0)^\top,& \sigma_4(x,y, z) &:= (0, \cos(2\pi x), 0 )^\top\\
        \sigma_5(x,y,z) &:= (\sin(2\pi z), 0, 0)^\top,&\ \sigma_6(x,y,z) &:= (\cos(2\pi z), 0, 0)^\top\\
        \sigma_7(x,y, z) &:= (0, 0, \sin(2\pi x))^\top,& \sigma_8(x,y, z) &:= (0, 0, \cos(2\pi x))^\top\\
        \sigma_9(x,y,z) &:= (0, \sin(2\pi z), 0)^\top,&\ \sigma_{10}(x,y,z) &:= (0, \cos(2\pi z), 0)^\top\\
        \sigma_{11}(x,y,z) &:= (0, 0, \sin(2\pi y))^\top,& \sigma_{12}(x,y,z) &:= (0, 0, \cos(2\pi y))^\top.
    \end{aligned}
\end{equation}

\begin{theorem}\label{thm:3d} The analogous statements of Theorem \ref{thm:lower_bound_low_modes} and its two corollaries hold for the three-dimensional model \eqref{spde:3D}.
\begin{enumerate}[(i)]
    \item     For $\kappa\geq 0$ and $f_0\in L_0^2(\tor^3) \setminus \{0\}$, let $f_t$ be the solution to \eqref{spde:3D}. The decay rate of the $L^2$-norm of the six innermost modes is lower bounded by
    \begin{equation}\label{ineq:LowerBound3D}
        \limsup_{t\to \infty} \frac{1}{t} \log \norm{\Pi_{\leq 1}f_t}_{L^2} \geq -(2\pi)^2 \left( 1 + \kappa \right ).
    \end{equation}
    \item The exponential dissipation rate $\lambda_{\omega, f_0}^\kappa$, as given by \eqref{eq:dissipation_rate}, is lower bounded by
    \begin{equation}
        \lambda_{\omega, f_0}^\kappa  \geq -(2\pi)^2 \left(1 + \kappa \right ).
    \end{equation}
    \item     For $s> 0$ the exponential mixing rate without diffusion is of order $1\land s$, that is, there exists a constant $C>0$ such that
    \begin{equation}
        C^{-1} (1\wedge s) \leq \gamma_s \leq C (1\wedge s), \quad \forall\, s>0.
    \end{equation}
\end{enumerate}
\end{theorem}
The proof of Theorem \ref{thm:3d} $(i)$ is almost analogous to the proof of Theorem \ref{thm:lower_bound_low_modes}. In Section \ref{sec:proof_3d}, we show how the proof of the 2-dimensional result can be adapted to establish Theorem \ref{thm:3d} $(i)$. Once $(i)$ is proven, $(ii)$ follows immediately. The proof of $(iii)$ is analogous to that of Corollary \ref{cor:mixingRate}. Like for the 4-modes model, \cite[Theorem 6.1]{gess2025stabilization} applies to the 12-modes model \eqref{spde:3D}. On the one hand, this proves a linear lower bound on $\gamma_s$ for small $s$, a necessary step for $(iii)$. On the other hand, their result establishes an upper bound $\limsup_{\kappa \to 0}\lambda^\kappa<0$. Together with $(ii)$, this shows that the Batchelor scale conjecture is verified for the three-dimensional model.

\subsection{Related literature}
The dependence on $\kappa$ of the almost-sure exponential dissipation rate $\lambda^\kappa$ has been studied from different angles throughout the literature. In \cite{miles2018diffusion} the dissipation rate has been estimated by means of numerical experiments. Their results indicate that $\lambda^\kappa$ is effectively independent of $\kappa$ for small $\kappa$, thus supporting the Batchelor scale conjecture. There are heuristic derivations of the Batchelor scale conjecture that so far did not yield rigorous proofs \cite{antonsen1996role, miles2018shell}.

The first uniform-in-diffusivity upper bound on the exponential dissipation rate has been established by Bedrossian, Blumenthal and Punshon-Smith \cite{bedrossian2021PTRF}. The authors show $\limsup_{\kappa \to 0} \lambda^\kappa<0$ for velocity fields given by the solutions to the 2-dimensional stochastic Navier--Stokes equation and by solutions to the 3-dimensional hyperviscous Navier--Stokes equation. Their proof uses ergodicity of the two-point motion of a random dynamical system corresponding to the advection-diffusion equation \eqref{pde:advdiff}. This technique goes back to Baxendale and Stroock  \cite{baxendale1988large} and has proven successful in establishing exponential mixing \cite{dolgopyat2004sample, bedrossian2022AOP}. Using similar methods, Gess and Yaroslavtsev \cite{gess2025stabilization} established a uniform upper bound on $\lambda^\kappa$ for white-in-time velocity fields, including the 4-modes model \eqref{spde:main2D} and its three-dimensional analogue \eqref{spde:3D}. 

\revision{For white-in-time models of the form \eqref{spde:advDiff}, Luo, Tang and Zhao \cite{luo2024elementary}, derive a quantitative upper bound for $\lambda_\kappa$ using elementary methods.} A different approach to proving a uniform-in-diffusivity upper bound on $\lambda^\kappa$ has been taken by Coti Zelati, Drivas and Gvalani \cite{zelati2023statistically} and Rowan \cite{rowan2024anomalous}. They consider white-in-time velocity fields and derive a formula for the expected $L^2$-norm of the scalar $f_t$ using its two-point correlation function. 

To the best of our knowledge, the only work that establishes a lower bound on the exponential dissipation rate $\lambda^\kappa$ is the recent work by Hairer, Punshon-Smith, Rosati and Yi \cite{hairer2024lower} where they establish $\lambda^\kappa \gtrsim \kappa^{-q}$ for any $q>3$ for velocity fields given by the 2-dimensional stochastic Navier--Stokes equation. Their proof is based on a concept they call \emph{high-frequency stochastic instability}.
For sufficiently noisy vector fields some Fourier mass of the solution to the advection-diffusion equation  is moved to low modes in a short amount of time with high probability. Once there is some mass in the low modes, diffusion kills the high modes quickly and decreases the filamentation length. This method heavily uses the fact that the noise forcing the velocity field is non-degenerate, that is, the noise acts on high modes. The 4-modes model \eqref{spde:main2D} is supported on the lowest four modes, such that high-frequency stochastic instability, as described in \cite{hairer2024lower}, cannot occur for this system. Our proof of the lower bound $\liminf_{\kappa \to 0} \lambda^\kappa>-\infty$ follows a different strategy. Instead of proving that Fourier mass is likely to return to low modes, we show that Fourier mass only ever leaves the low modes at an exponential rate that is uniform in $\kappa$; see Theorem \ref{thm:lower_bound_low_modes}. However, our proof heavily uses the fact that the 4-modes model is white-in-time and our method does not generalize to the 2-dimensional Navier--Stokes equation.

\revision{It is known in both discrete and continuous time that a generic compact linear cocycle, either has a collapsed Lyapunov spectrum, i.e.~the top Lyapunov exponent is $-\infty$, or it admits a \emph{dominated splitting} \cite{bessa2008lyapunov, bessa2019lyapunov}. This implies that any transfer operator cocycle of an advection-diffusion equation can be perturbed to either have collapsed Lyapunov spectrum, or to have a dominated splitting. We note that is unclear whether these perturbations can be chosen such that the perturbed cocycle still corresponds to the transfer operators of an advection-diffusion equation.}

\section{Proof of the Main Result}\label{sec:proof}
This section contains the proof of Theorem \ref{thm:lower_bound_low_modes}. In Section \ref{sec:proof_3d}, we show how the proof can be adapted to the three-dimensional model \eqref{spde:3D}, thus proving Theorem \ref{thm:3d}.

\revision{
In It\^o form, the general advection-diffusion SPDE \eqref{spde:advDiff} reads
\begin{equation}\label{spde:Ito}
    \dd f_t = \kappa \Delta f_t \, \dd t + \frac{1}{2} \sum_{j=1}^m (\sigma_j \cdot \nabla)^2 f_t \, \dd t - \sum_{j = 1}^m \sigma_j \cdot \nabla f_t \, \dd W_t^j.
\end{equation}
In the case of the 4-modes model \eqref{spde:main2D}, and the analogous 12-modes model \eqref{spde:3D}, the It\^o correction term equates to a multiple of $\Delta f_t$. This Laplace term does, however, not act as 'true diffusion' since for $\kappa=0$ the SPDE is a pure transport equation and preserves the $L^2$-norm of solutions. In Fourier coordinates, the advection term turns into a convolution with the Fourier coefficients of the velocity fields $\sigma_j$. For the 4-modes model, an elementary computation yields for $\mathbf{k}\in \Z^2$
\begin{equation}
    \dd \widehat{f}_\mathbf{k}(t) = -(2\pi)^2 \left(\frac{1}{2} + \kappa\right) \abs{\mathbf{k}}^2 \widehat{f}_\mathbf{k} (t) \, \dd t - 2\pi i \sum_{j=1}^4 \sum_{\mathbf{l} \in \Z^2} \big((\widehat{\sigma_j})_\mathbf{l} \cdot \mathbf{k}\big)\, \widehat{f}_{\mathbf{k}-\mathbf{l}}(t) \, \dd W_t^j.
\end{equation}
After inserting the Fourier coefficients of the four velocity fields \eqref{eq:4modes}, we obtain a closed representation of the evolution of the Fourier modes.
}

\begin{proposition}\label{prop:Fourier_representation}
    The evolution of the Fourier coefficients $\widehat{f}_{k,l}(t) \in \C$ in the 4-modes model \eqref{spde:main2D} is given by
    \begin{nalign}\label{eq:SDE_Fourier}
        \dd \widehat{f}_{k,l}(t) =  - (2 \pi)^2 \left(\frac{1}{2} + \kappa \right) (k^2 + l^2) \widehat{f}_{k,l}(t)\, \dd t -\pi  \Big( &k \big(\phantom{i} \widehat{f}_{k,l-1}(t) - \phantom{i}\widehat{f}_{k,l+1}(t)\big)\,\dd W_t^1 \\
        +\, &k \big(i \widehat{f}_{k,l-1}(t) + i\widehat{f}_{k,l+1}(t)\big)\,\dd W_t^2 \\
        +\, &\parbox{\widthof{k}}{l} \big(\phantom{i} \widehat{f}_{k-1,l}(t) - \phantom{i}\widehat{f}_{k+1,l}(t)\big)\,\dd W_t^3 \\
        + \,&\parbox{\widthof{k}}{l} \big(i \widehat{f}_{k-1,l}(t) + i\widehat{f}_{k+1,l}(t)\big)\,\dd W_t^4 \Big).
    \end{nalign}
\end{proposition}

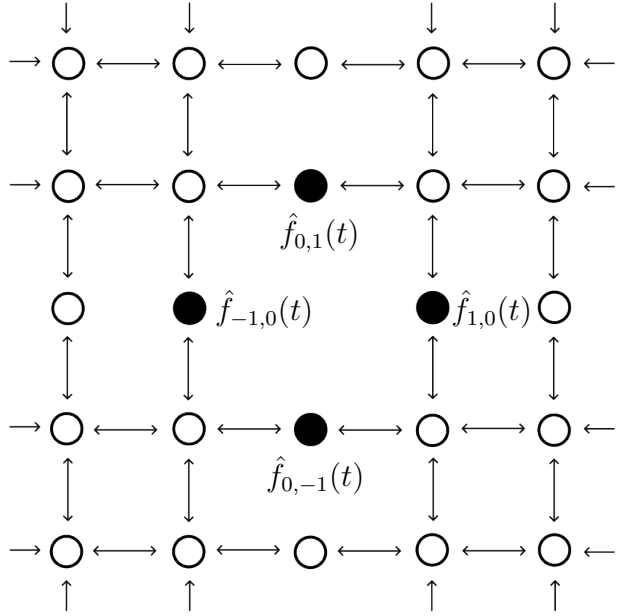
\begin{figure}
    \centering
    \begin{overpic}[width=0.5\linewidth]{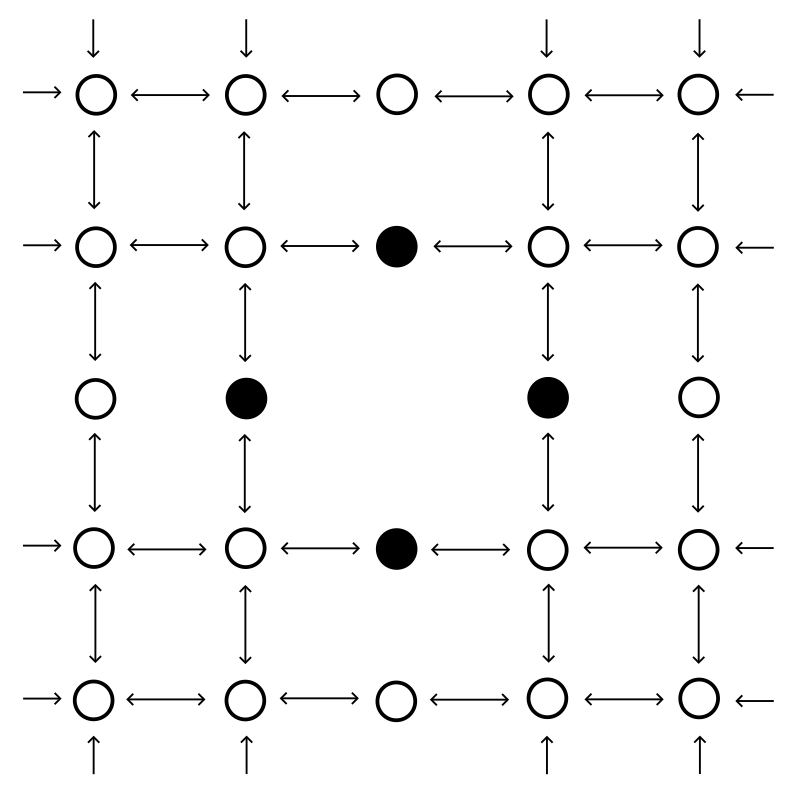}
        \put(45,60){$\hat f_{0,1}(t)$}
        \put(72,48){$\hat f_{1,0}(t)$}
        \put(35,48){$\hat f_{-1,0}(t)$}
        \put(43,22){$\hat f_{0,-1}(t)$}
    \end{overpic}
    \caption{Fourier modes $\hat f_{k,l}(t)$ of the solution $f_t$ to \eqref{spde:advDiff}. The inner modes with $k^2+l^2= 1$ are drawn as solid dots while all other modes are drawn as hollow dots. The arrows indicate which modes interact via the transport term; see Proposition \ref{prop:Fourier_representation}.}
    \label{fig:modes}
\end{figure}

\subsection{Proof of the Lower Bound} \label{se:proofLowerBound}
Let $f_0 \in L_0^2(\tor^2)$ with $\norm{f_0}_{L^2} = 1$. We first consider the case $\Pi_{\leq1} f_0 \neq 0$ and prove the almost-sure lower bound \eqref{ineq:MainLowerBound}. In Section \ref{sec:ProofIC} we show the same holds for all initial conditions $f_0\in L^2_0(\tor^2) \setminus \{0\}$.

Let $X_1, X_2, X_3, X_4, Y_1, Y_2, Y_3, Y_4$ be the real-valued stochastic processes for which
\revision{
\begin{align*}
    &\begin{aligned}
        \widehat f_{1,0}(t) &= X_1(t) + i X_2(t)\\
        \widehat f_{0,1}(t) &= X_3(t) + i X_4(t)
    \end{aligned}
    &&\text{and}&
    &\begin{aligned}
        \widehat f_{1,-1}(t) - \widehat f_{1,1}(t) &= Y_1(t) + i Y_2(t)\\
        \widehat f_{1,-1}(t) + \widehat f_{1,1}(t) &= Y_3(t) + i Y_4(t).
    \end{aligned}
\end{align*}
By the relation $\widehat{f}_{-k,-l}(t) = \overline{\widehat{f}_{k,l}(t)}$, we have
\begin{align*}
    &\begin{aligned}
        \widehat f_{-1,0}(t) &= X_1(t) - i X_2(t)\\
        \widehat f_{0,-1}(t) &= X_3(t) - i X_4(t)
    \end{aligned}
    &&\text{and}&
    &\begin{aligned}
        \widehat f_{-1,1}(t) - \widehat f_{1,1}(t) &= Y_1(t) - i Y_4(t)\\
        \widehat f_{-1,1}(t) + \widehat f_{1,1}(t) &= Y_3(t) - i Y_2(t).
    \end{aligned}
\end{align*}
By Proposition \ref{prop:Fourier_representation}, the vector-valued process $\mathbf{X}(t) \coloneqq (X_1(t),X_2(t),X_3(t),X_4(t))^\top$ satisfies the SDE
\begin{equation}\label{eq:A(t)_full}
    \dd \mathbf{X}(t) = - (2\pi)^2 \left(\frac{1}{2} + \kappa\right) \mathbf{X}(t)\,\dd t - \pi \begin{pmatrix}
    Y_1(t)&-Y_4(t)&0&0\\Y_2(t)&Y_3(t)&0&0\\0&0&Y_1(t)&Y_2(t)\\0&0&-Y_4(t)&Y_3(t)
\end{pmatrix}\dd \mathbf{W}(t),
\end{equation}
}
where $\mathbf{W}(t) \coloneqq (W_t^1,W_t^2, W_t^3, W_t^4)^\top$. Introducing the constant $\gamma \coloneqq (2\pi)^2 (\frac{1}{2} + \kappa)$ and denoting the matrix in \eqref{eq:A(t)_full} by $A(t)$, we shorten this SDE to
\begin{equation}
    \dd \mathbf{X}(t) = -\gamma \mathbf{X}(t) \,\dd t - \pi A(t) \,\dd \mathbf{W}(t).
\end{equation}

By the definition of $\mathbf{X}(t)$, we find $\norm{\Pi_{\leq 1} f_t}_{L^2}^2 = 2|\mathbf{X}(t)|^2$. Define
\begin{equation}\label{eq:def_lambda}
    \lambda_t \coloneqq \log \norm{\Pi_{\leq 1} f_t}_{L^2} = \frac{1}{2} \log\left( 2|\mathbf{X}(t)|^2 \right). 
\end{equation}
\revision{
\revision{Since $\norm{f_0}_{L^2}=1$ and the norm of $f_t$ can only decrease, $\lambda_t$ is non-positive for all $t\geq 0$.} Let us assume for now that that $\mathbf{X}(t)$ almost surely does not become $0$. Then, the quantity $\lambda_t$ is well-defined, and its evolution is given by It\^o's formula. We verify that $\mathbf{X}(t)$ does almost surely not hit $0$ below in step (I).} The first two derivatives of the function $F(\mathbf{X})\coloneqq \frac{1}{2} \log(2|\mathbf{X}|^2)$  are given by
\begin{nalign}\label{eq:Hessian_F}
    \nabla F(\mathbf{X}) &= |\mathbf{X}|^{-2} \mathbf{X}, \\
    \mathrm HF(\mathbf{X}) &= |\mathbf{X}|^{-4}\big(|\mathbf{X}|^2 \Id - 2 \mathbf{X} \mathbf{X}^\top\big).
\end{nalign}
Applying It\^o's formula yields
\begin{equation}\label{eq:Itos_forumla}
    \dd \lambda_t = - \gamma \nabla F(\mathbf{X}(t)) \cdot \mathbf{X}(t)\, \dd t + \underbrace{\frac{1}{2}\pi^2 \operatorname{tr} \left( A(t)^\top HF(\mathbf{X}(t)) A(t) \right)}_{\eqqcolon\mu(t)} \dd t - \pi \nabla F(\mathbf{X}(t))^\top A(t) \: \dd \mathbf{W}(t).
\end{equation}
\revision{Integrating in time, and using the fact that $\nabla F(\mathbf{X}(t)) \cdot \mathbf{X}(t) = 1$, we arrive at}
\begin{equation}\label{eq:lambda_integrated}
    \lambda_t = \lambda_0 -\gamma t + \int_0^t \mu(s) \:\dd s + \underbrace{(-\pi) \int_0^{t} \: |\mathbf{X}(s)|^{-2} \mathbf{X}(s)^\top A(s)\: \dd \mathbf{W}(s)}_{\eqqcolon M_t}.
\end{equation}
\revision{
The rest of the proof consists of three steps:
\begin{enumerate}[(I)]
    \item Show that the drift term $\mu(t)$ is non-negative.
    \item Show that $\mathbf{X}(t)$ almost surely does not become zero, such that \eqref{eq:lambda_integrated} is valid for all $t$.
    \item Show that the local martingale  satisfies $\limsup_{t\to \infty} \frac{1}{t} M_t \geq 0$.
\end{enumerate}
 Once these three steps are proven, \eqref{eq:lambda_integrated} yields
\begin{equation}
    \limsup_{t \to \infty} \frac{1}{t} \lambda_t \geq - \gamma.
\end{equation}
By definition of $\lambda_t$ in \eqref{eq:def_lambda}, and since $\gamma = (2\pi)^2 \big(\frac{1}{2} + \kappa\big)$, this is equivalent to
\begin{equation}
    \limsup_{t\to \infty} \frac{1}{t} \log \norm{ \Pi_{\leq 1} f_t}_{L^2} \geq -(2\pi)^2 \left(\frac{1}{2} + \kappa\right), 
\end{equation}
which is the desired bound \eqref{ineq:MainLowerBound}.\\
}

\textbf{(I) The drift term $\mu(t)$ is non-negative.} \\ 
We explicitly compute the drift term $\mu(t)$ as it is defined in \eqref{eq:Itos_forumla} using the formula of the Hessian $HF$ from \eqref{eq:Hessian_F}.
\begin{nalign}\label{eq:drift_comp}
    \mu(t) &= \frac{1}{2} \pi^2 |\mathbf{X}(t)|^{-2} \operatorname{tr} \left( A(t)^\top A(t) \right ) - \pi^2|\mathbf{X}(t)|^{-4} \operatorname{tr}\left(A(t)^\top \mathbf{X}(t) \mathbf{X}(t)^\top A(t) \right) \\
    &= \frac{1}{2} \pi^2 |\mathbf{X}(t)|^{-2} \left( \norm{A(t)}_{\text{Fr}}^2 - 2\left| A(t)^\top \frac{\mathbf{X}(t)}{|\mathbf{X}(t)|}\right|^2 \right) \\
    &\geq \frac{1}{2} \pi^2 |\mathbf{X}(t)|^{-2}\left( \norm{A(t)}_\text{Fr}^2 - 2 \norm{A(t)}_\text{op}^2 \right),
\end{nalign}
where $\norm{A(t)}_{\text{Fr}}= \operatorname{tr}(A(t)^\top A(t))^{1/2}$ is the Frobenius norm. This shows that a sufficient condition for the drift $\mu(t)$ to be non-negative is the algebraic condition
\revision{
\begin{equation} \label{eq:algebraic_condition}
    \norm{A(t)}_\text{Fr}^2 \geq 2 \norm{A(t)}_\text{op}^2.
\end{equation}
}
Observe that $A(t)$ is in block form with blocks $B_1, B_2 \in \R^{2\times 2}$; see \eqref{eq:A(t)_full}. Since $B_1$ and $B_2$ are transpose of one another, their Frobenius norms and operator norms coincide. Since the operator norm of a matrix is upper bounded by its Frobenius norm, we conclude
\begin{equation}
    \norm{A(t)}_\text{Fr}^2 = 2 \norm{B_i}_\text{Fr}^2 \geq 2 \norm{B_i}^2_\text{op} = 2 \norm{A(t)}_\text{op}^2.
\end{equation}
This verifies the condition \eqref{eq:algebraic_condition}, and shows that the drift term $\mu(t)$ is a non-negative, adapted process, no matter what the values $Y_i(t)$ are.\\

\textbf{(II) $\mathbf{X}(t)$ almost surely does not become zero.} \\
Define the stopping time
\begin{equation}
    \tau_n \coloneqq \inf \{t\geq 0 \mid \lambda_t \leq -n\}, \qquad \tau_\infty \coloneqq \lim_{n\to \infty} \tau_n.
\end{equation}
The time $\tau_\infty$ is the first time that $\mathbf{X}(t)$ becomes zero. Fix $n\in \N$ and stop the processes $\lambda_t$ and $\mathbf{X}(t)$ at time $\tau_n$. The stopped process $\lambda_{t\land \tau_n}$ is finite and we can apply It\^o's formula. Hence, \eqref{eq:lambda_integrated} is valid up until time $\tau_n$. The stopped process $\lambda_{t\land \tau_n}$ is given by
\begin{equation}\label{eq:lambda_integrated_tau}
    \lambda_{t\land \tau_n} = \lambda_0 -\gamma (t\land \tau_n) + \int_0^{t\land \tau_n} \mu(s) \:\dd s + M_{t\land \tau_n}.
\end{equation}
In the previous step, we proved that $\mu(s)\geq 0$. Thus, the expected value of $\lambda_{t \land \tau_n}$ is bounded from below by
\begin{equation}
    \E[\lambda_{t \land \tau_n}]  \geq \lambda_0 - \gamma t. 
\end{equation}
From Markov's inequality for $-\lambda_{t \land \tau_n}$, we get
\begin{equation}
- \lambda_0 + \gamma t \geq \E[-\lambda_{t \land \tau_n}] \geq \P(\tau_n \leq t)\, n.
\end{equation}
For fixed $t\geq 0$, we let $n\to \infty$. This shows that $\tau_\infty > t$ with probability $1$, and $\mathbf{X}(t)$ almost surely never becomes zero. \\

\textbf{(III) It holds $\limsup_{t\to \infty} \frac{1}{t} M_t \geq 0$.}\\
Consider the local martingale term $M_t$ in \eqref{eq:lambda_integrated}. Its quadratic variation is given by
\begin{equation}
    \langle M \rangle_t = \pi^2 \int_0^t \: |\mathbf{X}(s)|^{-4} |\mathbf{X}(s)^\top A(s)|^2 \: \dd s.
\end{equation}
The total quadratic variation is given by $\langle M \rangle_\infty= \lim_{t\to \infty} \langle M \rangle_t$. If $\langle M \rangle_\infty < \infty$, the limit $M_\infty = \lim_{t\to \infty} M_t$ exists. For $t\geq 0$, define the random time-change
\begin{equation}
    T_t \coloneqq \inf \{T\geq 0 \mid \langle M \rangle_T \geq t \}.
\end{equation}
As long as $t < \langle M \rangle_\infty$, the time $T_t$ is finite. We construct the process
\begin{equation}
    B_t \coloneqq \begin{cases}
        M_{T_t} \quad &\text{if} \quad t < \langle M \rangle_\infty, \\
        M_\infty\quad &\text{if} \quad t \geq \langle M \rangle_\infty.
    \end{cases}
\end{equation}
We use a generalization of the Dambis, Dubins--Schwarz Theorem which is formulated in \cite[Theorem V.1.7]{rezuv1999continuous}. The theorem asserts that $B_t$ is a Brownian motion that is stopped at $\langle M \rangle_\infty$. If $W_t$ is a Brownian motion with respect to the filtration $\{\mathcal{F}_t\}_{t\geq 0}$, then $B_t$ is a stopped Brownian motion with respect to the filtration $\{F_{T_t}\}_{t\geq 0}$. A stopped Brownian motion either remains bounded (in case it stops) or becomes positive infinitely often. In the first case, which occurs when $\langle M\rangle_\infty < \infty$, the process $M_t$ stays bounded and we have
\begin{equation}
    \lim_{t\to \infty} \frac{1}{t} M_t = 0.
\end{equation}
In the latter case, $B_t$ becomes positive infinitely often and we find
\begin{equation}
    \limsup_{t\to \infty} \frac{1}{t} M_t \geq \limsup_{t\to \infty} \frac{1}{T_t} M_{T_t} = \limsup_{t\to \infty} \frac{1}{T_t} B_t \geq 0.
\end{equation}

This proves the last of the three steps, and the proof is complete.

\subsection{Independence of the initial condition} \label{sec:ProofIC}
In the previous section, we proved that for initial conditions $f_0 \in L^2_0(\tor^2)$ with $\Pi_{\leq 1} f_0 \neq 0$, we almost surely have
\begin{equation}
    \limsup_{t \to \infty} \frac{1}{t} \log \norm{\Pi_{\leq 1} f_t}_{L^2} \geq - (2\pi)^2 \left( \frac{1}{2} + \kappa \right ).
\end{equation}
In this section, we show that for all $f_0 \in L^2_0(\tor^2) \setminus \{0\}$ the same almost-sure lower bound holds. We need the following lemma.

\begin{lemma}\label{lem:initial_condition}
    Let $f_0 \in L^2_0(\tor^2) \setminus \{0\}$. For every $(k,l)\in \Z^2 \setminus \{0\}$ we have 
    \begin{equation}
        \P\big(\widehat{f}_{k,l}(t) = 0, \forall\, t\geq 0\big) = 0.
    \end{equation}
\end{lemma}

\begin{proof}
    Fix $f_0 \in L^2_0(\tor^2) \setminus \{0\}$. Consider the set of modes which have a positive probability to be zero for all times,
    \begin{equation}
        N \coloneqq \left\{(k,l) \in\Z^2 \setminus \{0\} \: \big| \: \P\big(\widehat{f}_{k,l}(t) = 0, \forall \, t\geq 0\big) > 0\right \}.
    \end{equation}
    
    Recall the evolution of the Fourier modes $\widehat{f}_{k,l}(t)$; see~\eqref{eq:SDE_Fourier}. 
    Let $V_{k,l}$ denote the sum of the quadratic variations of the real and imaginary parts of $\hat f_{k,l}(t)$, which is given by
    \begin{nalign}\label{eq:quadratic_variation}
        V_{k,l}(t) &\coloneqq\big \langle \Re(\widehat{f}_{k,l})\big \rangle_t +\big \langle \Im(\widehat{f}_{k,l})\big \rangle_t \\
        &\,= 2\pi^2\int_0^t \left[k^2\big(|\widehat{f}_{k,l-1}(s)|^2 + |\widehat{f}_{k,l+1}(s)|^2 \big)  + l^2\big(|\widehat{f}_{k-1,l}(s)|^2 + |\widehat{f}_{k+1,l}(s)|^2 \big) \right]\,\dd s.
    \end{nalign}
    We say two modes $(k_1, l_1), (k_2, l_2) \in \Z^2 \setminus \{0\}$ are \emph{adjacent} if
    \begin{equation}
       \big(k_1 =k_2 \neq 0,
        \quad   l_1 = l_2 \pm 1\big) \qquad \text{ or } \qquad \big( k_1 = k_2 \pm 1,\quad   l_1 = l_2 \neq 0 \big);
    \end{equation}
    see Figure \ref{fig:modes}. If $(k,l) \in N$, then $\mathbb P(V_{k,l}(t) =0, \forall\, t\geq 0) >0$, and by \eqref{eq:quadratic_variation}, we have $(k',l')\in N$ for every $(k',l')$ adjacent to $(k,l)$. Since the graph with vertex set $\mathbb Z^2\setminus \{0\}$ and edges between all adjacent modes is connected, we must have either $N = \emptyset$ or $N = \mathbb Z^2\setminus\{0\}$. By assumption $N\neq \mathbb Z^2\setminus \{0\}$, so $N = \emptyset$.
\end{proof}

We are now able to prove Theorem \ref{thm:lower_bound_low_modes}.

\begin{proof}[Proof of Theorem \ref{thm:lower_bound_low_modes}]
    Let $f_0 \in L^2_0(\tor^2) \setminus \{0\}$. By Lemma \ref{lem:initial_condition}, almost surely, there is some time $t\geq 0$ such that $\Pi_{\leq 1}f_t \neq 0$. For $\epsilon>0$, define the stopping time
\begin{equation}\label{eq:tau_epsilon}
    \tau_\epsilon \coloneqq\inf \{t>0 \mid\norm{\Pi_{\leq 1}f_t}_{L^2} > \epsilon\}.
\end{equation}
If $\tau_\epsilon<\infty$, we have $\Pi_{\leq 1} f_{\tau_\epsilon} \neq 0$. Using the result of Section \ref{se:proofLowerBound} and the strong Markov property, we get 
\begin{equation}
    \mathbb P\left(\limsup_{t \to \infty} \frac{1}{t} \log \norm{\Pi_{\leq 1}f_t}_{L^2} \geq -(2\pi)^2\left(\frac{1}{2} + \kappa\right ) \Bigg|\: \tau_\epsilon<\infty \right) = 1,
\end{equation}
and thus
\begin{equation}
    \P\left(\limsup_{t \to \infty} \frac{1}{t} \log \norm{\Pi_{\leq 1}f_t}_{L^2} \geq -(2\pi)^2\left(\frac{1}{2} + \kappa\right)\right) \geq \P(\tau_\epsilon < \infty).
\end{equation}
Note that the left-hand side does not depend on $\epsilon$. When taking $\epsilon\to 0$ the right-hand side converges to 1 since $\Pi_{\leq 1} f_t$ is almost surely non-zero for some $t\geq 0$. Thus \eqref{ineq:MainLowerBound} holds almost surely.

\end{proof}

\subsection{Three-dimensional model}\label{sec:proof_3d}
This section is devoted to the proof of Theorem \ref{thm:3d}. The proof is almost analogous to that of the two-dimensional model. Therefore, we omit the parts of the proof that are the same and elaborate which steps need to be adapted. As before, we first rewrite the SPDE \eqref{spde:3D} in Fourier coordinates. Note that on the three-dimensional torus, the Fourier modes are indexed by three integers $k,l,m \in \Z$.
\begin{proposition}\label{prop:3d_Fourier}
    The evolution of the Fourier coefficients $\widehat{f}_{k,l,m}(t) \in \C$ of the 12-modes model \eqref{spde:3D} is given by\revision{
    \begin{nalign}
        \dd \widehat{f}_{k,l,m}(t) =  - (2 \pi&)^2 \left(1+\kappa \right) (k^2 + l^2 + m^2) \widehat{f}_{k,l,m}(t) \, \dd t  \\
        -\pi \Big( \:\phantom{+}&\parbox{\widthof{m}}{k} \big( \widehat{f}_{k,l-1, m}(t) - \widehat{f}_{k,l+1, m}(t)\big)\, \dd W_t^1 \hspace{5pt} + \,\parbox{\widthof{m}}{k} \big( i \widehat{f}_{k,l-1, m}(t) + i\widehat{f}_{k,l+1, m}(t)\big)\, \dd W_t^2 \\
        +&\parbox{\widthof{m}}{l} \big( \widehat{f}_{k-1,l, m}(t) - \widehat{f}_{k+1,l, m}(t)\big)\, \dd W_t^3\hspace{5pt}  + \,\parbox{\widthof{m}}{l} \big( i \widehat{f}_{k-1,l, m}(t) + i\widehat{f}_{k+1,l, m}(t)\big)\, \dd W_t^4 \\
        +&\parbox{\widthof{m}}{k} \big( \widehat{f}_{k,l, m-1}(t) - \widehat{f}_{k,l, m+1}(t)\big)\, \dd W_t^5\hspace{5pt}  + \,\parbox{\widthof{m}}{k} \big( i \widehat{f}_{k,l, m-1}(t) + i\widehat{f}_{k,l, m+1}(t)\big)\, \dd W_t^6 \\
        +&\parbox{\widthof{m}}{m} \big( \widehat{f}_{k-1,l, m}(t) - \widehat{f}_{k+1,l, m}(t)\big)\, \dd W_t^7\hspace{5pt}  + \,\parbox{\widthof{m}}{m} \big( i \widehat{f}_{k-1,l, m}(t) + i\widehat{f}_{k+1,l, m}(t)\big)\, \dd W_t^8 \\
        +&\parbox{\widthof{m}}{l} \big( \widehat{f}_{k,l, m-1}(t) - \widehat{f}_{k,l, m+1}(t)\big)\, \dd W_t^9\hspace{5pt} +\, \parbox{\widthof{m}}{l} \big( i \widehat{f}_{k,l, m-1}(t) + i\widehat{f}_{k,l, m+1}(t)\big)\, \dd W_t^{10} \\
        +&m \big( \widehat{f}_{k,l-1, m}(t) - \widehat{f}_{k,l+1, m}(t)\big)\, \dd W_t^{11} +\, m \big( i \widehat{f}_{k,l-1, m}(t) + i\widehat{f}_{k,l+1, m}(t)\big)\, \dd W_t^{12} \: \Big). \\
    \end{nalign}}
\end{proposition}

Assume that $\Pi_{\leq 1} f_0 \neq 0$, and let $X_1, \hdots, X_6 \in \R$ be the real-valued stochastic processes for which
\begin{align}
    \widehat{f}_{1,0,0} (t) &= X_1(t) + iX_2(t) \\
    \widehat{f}_{0,1,0} (t) &= X_3(t) + iX_4(t) \\
    \widehat{f}_{0,0,1} (t) &= X_5(t) + iX_6(t).
\end{align}
Based on Proposition \ref{prop:3d_Fourier} the evolution of the vector $\mathbf{X}(t) \coloneqq (X_1(t) ,\hdots, X_6(t))^\top$ can be written as
\begin{equation}
    \dd \mathbf{X}(t) = -\gamma \mathbf{X}(t) \: \dd t - \pi A(t) \: \dd \mathbf{W}(t),
\end{equation}
where $\gamma \coloneqq \revision{(2\pi)^2(1+\kappa)}$,  $\mathbf{W}(t) = (W_t^1, \hdots W_t^{12} )^\top$, and $A(t) \in \R^{6 \times 12}$. The crucial step in the proof of the two-dimensional model was to verify the algebraic condition \eqref{eq:algebraic_condition}, which requires the squared Frobenius norm of $A(t)$ to be at least twice its squared operator norm. We verify this condition for the three-dimensional model.

To write the matrix $A(t)$, define the variables $Y_1, \hdots, Y_{12} \in \R$ such that\revision{
\begin{align*}
    \widehat{f}_{1,-1,0}(t) - \widehat{f}_{1,1,0}(t) &= Y_1(t) + i Y_2(t),  & \widehat{f}_{1,-1,0}(t) + \widehat{f}_{1,1,0}(t) &= Y_3(t) + i Y_4(t) \\
    \widehat{f}_{1,0,-1}(t) - \widehat{f}_{1,0,1}(t) &= Y_5(t) + i Y_6(t),  &  \widehat{f}_{1,0,-1}(t) + \widehat{f}_{1,0,1}(t) &= Y_7(t) + i Y_8(t) \\
    \widehat{f}_{0,1,-1}(t) - \widehat{f}_{0,1,1}(t) &= Y_{9}(t) + i Y_{10}(t),  &  \widehat{f}_{0,1,-1}(t) + \widehat{f}_{0,1,1}(t) &= Y_{11}(t) + i Y_{12}(t).
\end{align*}
The function $f_t$ is real so that $\widehat{f}_{-k,-l,-m}(t) = \overline{\widehat{f}_{k,l,m}(t)}$. The matrix $A(t)$ is given by

\begin{equation}\label{eq:big_matrix}
    A(t) =
    \left(
    \begin{array}{cccc|cccc|cccc}
    Y_1 & -Y_4 & 0 & 0 & Y_5 & -Y_8 & 0 & 0 & 0 & 0 & 0 & 0 \\[2pt]
    Y_2 & Y_3 & 0 & 0 & Y_6 & Y_7 & 0 & 0 & 0 & 0 & 0 & 0 \\[2pt]
    0 & 0 & Y_1 & Y_2 & 0 & 0 & 0 & 0 & Y_{9} & -Y_{12} & 0 & 0 \\[2pt]
    0 & 0 & -Y_4 & Y_3 & 0 & 0 & 0 & 0 & Y_{10} & Y_{11} & 0 & 0 \\[2pt]
    0 & 0 & 0 & 0 & 0 & 0 & Y_5 & Y_6 & 0 & 0 & Y_{9} & Y_{10} \\[2pt]
    0 & 0 & 0 & 0 & 0 & 0 & -Y_8 & Y_7 & 0 & 0 & -Y_{12} & Y_{11}
    \end{array}
    \right),
\end{equation}
}
where we omitted the time-dependence of the $Y_i$. The structure of \eqref{eq:big_matrix} is similar to three copies of the matrix \eqref{eq:A(t)_full} from the two-dimensional example. 
We consider the matrix in block form $A(t) = (B \: |\: C \: | \: D\:)$, where each block is a ${6 \times 4}$ matrix.
By the argument from the two-dimensional proof, the squared Frobenius norm of each of these blocks is at least twice its squared operator norm. We conclude
\begin{equation}
    \norm{A(t)}_\text{Fr}^2 = \norm{B}_\text{Fr}^2 +  \norm{C}_\text{Fr}^2 +  \norm{D}_\text{Fr}^2 \geq 2 \norm{B}_\text{op}^2 +  2\norm{C}_\text{op}^2 +  2\norm{D}_\text{op}^2 \geq 2 \norm{A(t)}_\text{op}^2.
\end{equation}
From here on, the proof is analogous to the two-dimensional case. We obtain the almost-sure lower bound
\begin{equation*}
    \limsup_{t\to \infty} \frac{1}{t} \log \revision{\norm{\Pi_{\leq 1} f_t}_{L^2}} \geq -\gamma = -(2\pi)^2(1 + \kappa).
\end{equation*}
Showing that this lower-bound holds for all initial conditions $f_0 \in L^2(\tor^3)$ is analogous to the proof in Section \ref{sec:ProofIC}.

\appendix

\section{Estimates on the mixing rate $\gamma_s$ of the transport SPDE}\label{app:mixing rates}
In this appendix, we show estimates for the mixing rate $\gamma_s$ of the stochastic transport equation
\begin{equation}\label{spde:Transport2}
    \dd f_t = - \sum_{j = 1}^m \sigma_j \cdot \nabla f_t \circ \dd W_t^j.
\end{equation}
For $m \in \mathbb N$ and some $\alpha \in (0,1]$, $\sigma_1,\dots, \sigma_m \in C^{2, \alpha}(\mathbb T^d, \mathbb R^d)$ are autonomous, divergence-free vector fields on the torus, and $W^1, \dots, W^m$ are independent Brownian motions on the stochastic basis $(\Omega, \mathcal F, \mathcal F_t, \mathbb P)$.  Weak-in-space, strong-in-probability solutions can be constructed using a Lagrangian approach. The Lagrangian SDE corresponding to \eqref{spde:Transport2} is given by
\begin{equation}\label{sde:Lagrangian}
    \dd \mathbf X_t = \sum_{j = 1}^m \sigma_j(\mathbf X_t) \circ \dd W_t^j.
\end{equation}
By a classical result due to Kunita (see, e.g., \cite{kunita2006stochastic}), there exists a stochastic flow of diffeomorphism, that is, a map $\varphi: (\omega, t_0,t,\mathbf x) \mapsto \varphi^{t_0,t}_\omega(\mathbf x)$, defined for $\omega \in \Omega, 0\leq t_0 \leq t < \infty$ and $\mathbf x \in \mathbb T^d$ with the following properties.
\begin{enumerate}[(i)]
    \item For every $t_0\geq 0$ and every $\mathbf x \in \mathbb T^d$, the stochastic process $\mathbf X_t := \varphi^{t_0,t}_\omega(\mathbf x), t \geq t_0$ is a strong solution to \eqref{sde:Lagrangian} with $\mathbf X_t =\mathbf x$.
    \item For almost every $\omega \in \Omega$, the map $\varphi^{t_0,t}_\omega: \mathbb T^d \to \mathbb T^d$ is a diffeomorphism for all $0 \leq t_0 \leq t$. 
    \item For almost every $\omega \in \Omega$ and for every $0 \leq t_0 \leq t_1 \leq t_2$, we have $\varphi_\omega^{t_0, t_2} = \varphi_\omega^{t_1, t_2} \circ \varphi_\omega^{t_0, t_1}$ and $\varphi^{t_0, t_0}_\omega = \operatorname{Id}_{\mathbb T^d}$.
\end{enumerate}

Since the vector fields $\sigma_1, \dots, \sigma_m$ are divergence-free, the stochastic flow $\varphi$ is volume preserving, that is, almost surely $\det \mathrm D_{\mathbf x} \varphi^{t_0,t}_\omega (\mathbf x) = 1$ for all $t_0 \leq t$ and $\mathbf x \in \tor^d$. The solution $f_t \in L^{\infty}([0,\infty), L^2(\mathbb T^d, \mathbb R))$ to \eqref{spde:Transport2} with initial state $f_0 \in L^2(\mathbb T^d, \mathbb R)$ can now be constructed by 
\begin{equation}\label{eq:ftConstruction}
    f_t = f_0 \circ (\varphi^{0,t}_\omega)^{-1}.
\end{equation}
For any test function $g \in C^\infty(\mathbb T^d, \mathbb R)$ we have
$$g(\mathbf X_t) = g(\mathbf x) + \sum_{j = 1}^m \int_0^t \sigma_j(\mathbf X_r) \cdot \nabla g(\mathbf X_r) \circ \dd W^j_{r},$$
where $\mathbf X_t := \varphi^{0,t}_\omega(\mathbf x)$. Thus, we can compute
\begin{align*}
    \langle f_t, g\rangle_{L^2} &= \left\langle  f_0 \circ (\varphi^{0,t}_\omega)^{-1}, g\right\rangle_{L^2} = \left\langle  f_0 , g\circ \varphi^{0,t}_\omega\right\rangle_{L^2}\\
    &= \left\langle  f_0 , g\right\rangle_{L^2} + \sum_{j = 1}^m \int_0^t \left\langle  f_0 , (\sigma_j \cdot \nabla g)\circ \varphi^{0,r}_\omega \right\rangle_{L^2} \circ \dd W_r^j\\
    &= \left\langle  f_0 , g\right\rangle_{L^2} + \sum_{j = 1}^m \int_0^t \left\langle  f_r ,  \nabla \cdot(g \sigma_j)\right\rangle_{L^2} \circ \dd W_r^j,
\end{align*}
that is, $f_t$ is indeed a weak-in-space, strong-in-probability solution to \eqref{spde:Transport2}.

We use \eqref{eq:ftConstruction} to derive an upper bound for the mixing rate $\gamma_s$. By \cite[Proposition 2.1]{baxendale1989lyapunov}, the stochastic flow of diffeomorphism $\varphi$ satisfies 
$$\mathbb E\left[\sup_{0\leq t \leq 1, \mathbf x \in \mathbb T^d}\left\|\mathrm D_{\mathbf x}(\varphi^{0,t}_\omega)^{-1}(\mathbf x)\right\|\right] < \infty.$$
Consequently, by Kingman's subadditive ergodic theorem \cite{Kingman}, there exists a finite deterministic number $\Lambda \geq 0$ , such that
$$\Lambda = \lim_{t \to \infty} \frac{1}{t} \sup_{\mathbf x\in \mathbb T^d}\log \|(\mathrm D_{\mathbf x}\varphi^{0,t}_\omega(\mathbf x))^{-1}\|, \quad \text{almost surely.}$$
\begin{proposition}\label{prop:GammaSUpperBound}
    For all $s>0$, if \eqref{spde:Transport2} is exponentially mixing for $s$, the mixing rate is upper bounded by $\gamma_s \leq \Lambda s$.
\end{proposition}
\begin{proof}
The growth of the $H^1$-norm of an arbitrary initial state $f_0 \in H^1$ can be controlled by
    \begin{align*}
        \|f_t\|_{H^1} &= \left(\int_{\mathbb T^d} |\nabla f_t(\mathbf x)|^2\, \dd \mathbf x\right)^{\frac 12} = \left(\int_{\mathbb T^d} \left|\nabla \left(f_0\circ (\varphi^{0,t}_\omega)^{-1}\right)(\mathbf x)\right|^2\, \dd \mathbf x\right)^{\frac 12}\\
        &\leq \left(\sup_{\mathbf x\in \mathbb T^d}\left\|(\mathrm D_{\mathbf x}\varphi_\omega^{0,t}(\mathbf x))^{-1}\right\|\right) \left(\int_{\mathbb T^d}  \left|\nabla f_0\left( (\varphi^{0,t}_\omega)^{-1}(\mathbf x)\right)\right|^2\, \dd \mathbf x\right)^{\frac 12}\\
        &\leq \tilde C_\omega\, e^{(\Lambda+\epsilon) t}\,\|f_0\|_{H^1},
    \end{align*}
    where $ \tilde C_\omega$ is a random constant given by 
    $$ \tilde C_\omega := \sup_{t'\geq 0} e^{-(\Lambda +\epsilon) t'}\sup_{\mathbf x\in \mathbb T^d}\left\|(\mathrm D_{\mathbf x}\varphi_\omega^{0, t'}(\mathbf x))^{-1}\right\|<\infty,$$
    which is finite by the definition of $\Lambda$. For any $s>0$, rearranging the interpolation inequality (Lemma \ref{lemm:interpolation} (i)) ${\|\cdot\|_{L^2} \leq \|\cdot\|_{H^{-s}}^{1/(1+s)}\|\cdot\|_{H^1}^{s/(1+s)}}$ gives $\|\cdot\|_{H^{-s}} \geq \|\cdot\|_{L^2}^{1+s} \|\cdot\|_{H^1}^{-s}$. For any $f_0 \in H^{1}$, we get
    $$\|f_t\|_{H^{-s}} \geq \|f_t\|_{L^2}^{1+s} \|f_t\|_{H^1}^{-s} \geq  \tilde C_\omega^{-s}\, e^{-s(\Lambda+\epsilon) t}\, \|f_0\|_{L^2}^{1+s}\,\|f_0\|^{-s}_{H^1}.$$
    In particular, the inequality \eqref{ineq:random_mixing} can not be satisfied for $\gamma = s\Lambda +2s\epsilon$, and thus $\gamma_s \leq s\Lambda +2s\epsilon$. Taking $\epsilon \to 0$ gives $\gamma_s \leq s\Lambda$, as claimed.
\end{proof}

We conclude the appendix showing that \eqref{spde:Transport2} is mixing for all $s>0$ given that it is mixing for some $s_0>0$. Furthermore we provide lower bounds for the mixing rates $\gamma_s$ in terms of the mixing rate $\gamma_{s_0}$. The following proposition is formulated in the context of the transport SPDE \eqref{spde:Transport2}, but the proof does not use the form of \eqref{spde:Transport2}. In particular, the statement is also true for random transport PDEs of the form \eqref{pde:Transport} and the proofs are identical.
\begin{proposition}[Lemma 5.1 in \cite{dewitt2025equivalence}] \label{prop:GammaSLowerBound}
    Let $s_0>0$ and suppose that \eqref{spde:Transport2} is exponentially mixing in the sense of \eqref{ineq:random_mixing} with mixing rate $\gamma_{s_0}$. Then \eqref{spde:Transport2} mixes exponentially for all $s>0$ and we have the following lower bounds for the mixing rate.
    \begin{enumerate}[(i)]
        \item $\gamma_s \geq \gamma_{s_0}$, for all $s>s_0$.
        \item $\gamma_s \geq \frac{s}{2s_0-s}\gamma_{s_0}$, for all $0<s<s_0$. 
    \end{enumerate}
\end{proposition}
The proof relies heavily on the following standard interpolation theorem. 
\begin{lemma}\label{lemm:interpolation}
Let $-\infty < s_1< s_2< s_3< \infty$.
\begin{enumerate}[(i)]
    \item For every $f \in H^{s_3}_0$, we have
        $$\|f\|_{H^{s_2}} \leq \|f\|_{H^{s_1}}^{(s_3-s_2)/(s_3-s_1)}\,\|f\|_{H^{s_3}}^{{(s_2-s_1)/(s_3-s_1)}}.$$
    \item For every $f \in H^{s_2}_0$ and every $\epsilon>0$, there exists an $f^\epsilon \in H^{s_3}_0$, such that
    \begin{equation}\label{ineqs:approx}
       \|f-f^\epsilon\|_{H^{s_1}} \leq \epsilon \|f\|_{H^{s_2}}, \quad  \|f^\epsilon\|_{H^{s_2}} \leq \|f\|_{H^{s_2}}, \quad \|f^\epsilon\|_{H^{s_3}} \leq \epsilon^{-(s_3-s_2)/(s_2-s_1)} \|f\|_{H^{s_2}}.
    \end{equation}
\end{enumerate} 
\end{lemma}
The fact that these inequalities hold without the appearance of constants is due to our choice of the $H^s$-norms \eqref{eq:Sobolev_norm} among all equivalent norms. We sketch the proofs.
\begin{proof}[Sketch of Proof]
    (i) is a direct consequence of \eqref{eq:Sobolev_norm} and the Hölder inequality. For (ii), let $R := \epsilon^{-1/(s_2-s_1)}$ and let
    $$f^\epsilon := \Pi_{\leq R} f := \sum_{\mathbf k \in \mathbb Z^d\setminus \{0\}, |\mathbf k| \leq R} \hat f_{\mathbf k} F_{\mathbf k}.$$
    Then 
    $$f-f^\epsilon = \Pi_{> R} f := \sum_{\mathbf k \in \mathbb Z^d\setminus \{0\}, |\mathbf k| > R} \hat f_{\mathbf k} F_{\mathbf k}$$
    and the inequalities in \eqref{ineqs:approx} follow from \eqref{eq:Sobolev_norm} and an elementary computation.
\end{proof}

\begin{proof}[Proof of Proposition \ref{prop:GammaSLowerBound}]
    The case $s>s_0$ is a trivial consequence of the embeddings $H^s \hookrightarrow H^{s_0}$ and $H^{-s_0} \hookrightarrow H^{-s}$. Let $0<s<s_0$ and $\gamma <\gamma_{s_0}$. By the definition of the mixing rate $\gamma_{s_0}$, there exists a random constant $C_\omega$, such that almost surely
    \begin{equation} \label{ineq:mixingAppendix}
        \|f_t\|_{H^{-s_0}} \leq C_{\omega} e^{-\gamma t} \|f_0\|_{H^{s_0}},\quad\forall\, f_0 \in H^{s_0}_0.
    \end{equation}
    We show that there exists a random constant $C'_{\omega}$ such that
    $$\|f_t\|_{H^{-s}} \leq C'_{\omega} e^{-s\gamma t/(2s_0-s)} \|f_0\|_{H^s},\quad\forall\, f_0 \in H^s_0.$$
    This implies $\gamma_s > \frac{s}{(2s_0-s)} \gamma$, and taking $\gamma \uparrow \gamma_{s_0}$ yields the claimed bound.
    
    Fix $f_0 \in H_0^s$ and $t>0$. By Lemma \ref{lemm:interpolation} (ii), for any $\epsilon > 0$, there exists a function $f_0^\epsilon \in H^{s_0}_0$ with
    $$\|f_0-f_0^\epsilon\|_{L^2} \leq \epsilon \|f_0\|_{H^{s}}, \quad  \|f_0^\epsilon\|_{L^2} \leq \|f_0\|_{L^2}, \quad \text{and} \quad \|f_0^\epsilon\|_{H^{s_0}} \leq \epsilon^{-(s_0-s)/s} \|f_0\|_{H^{s}}.$$
   Let $f_t^\epsilon$ denote the solution to \eqref{spde:Transport2} with initial state $f_0^\epsilon$. By \eqref{ineq:mixingAppendix}, we have 
    \begin{equation}\label{ineq:AppendixBound1}
        \|f^\epsilon_t\|_{H^{-s_0}} \leq C_{\omega} e^{-\gamma t} \|f^\epsilon_0\|_{H^{s_0}} \leq  C_{\omega} e^{-\gamma t} \epsilon^{-(s_0-s)/s} \|f_0\|_{H^{s}}.
    \end{equation}
    Using the fact that \eqref{spde:Transport2} preserves the $L^2$-norm, we furthermore get
    \begin{equation}\label{ineq:AppendixBound2}
        \|f^\epsilon_t\|_{L^2} = \|f^\epsilon_0\|_{L^2} \leq  \|f_0\|_{L^2}\leq \|f_0\|_{H^s}.
    \end{equation}
 The interpolation between the bounds \eqref{ineq:AppendixBound1} and \eqref{ineq:AppendixBound2}, using Lemma \ref{lemm:interpolation} (i), gives 
    \begin{equation} \label{ineq:feps}
        \|f^\epsilon_t\|_{H^{-s}} \leq \|f_t^\epsilon\|_{H^{-s_0}}^{s/s_0} \|f_t^\epsilon\|_{L^2}^{1-s/s_0} \leq C_{\omega}^{s/s_0} e^{-(s/s_0)\gamma t} \epsilon^{-(s_0-s)/s_0} \|f_0\|_{H^{s}}.
    \end{equation}
    Since \eqref{spde:Transport2} is linear, $f_t -f_t^\epsilon$ is also a solution with initial condition $f_0 - f_0^\epsilon$ and its $L^2$-norm is constant. This allows us to bound
    \begin{equation}\label{ineq:f-feps}
         \|f_t - f^\epsilon_t\|_{H^{-s}} \leq \|f_t - f^\epsilon_t\|_{L^2} =  \|f_0 - f^\epsilon_0\|_{L^2} \leq \epsilon \|f_0\|_{H^{s}}.
    \end{equation}
    Combining the estimates \eqref{ineq:feps} and \eqref{ineq:f-feps}, we get
    $$\|f_t\|_{H^{-s}} \leq  \|f^\epsilon_t\|_{H^{-s}} + \|f_t - f^\epsilon_t\|_{H^{-s}} \leq \left(C_{\omega}^{s/s_0} e^{-(s/s_0)\gamma t} \epsilon^{-(s_0-s)/s_0} + \epsilon\right) \|f_0\|_{H^{s}}.$$
    Choosing $\epsilon := e^{- s\gamma t/(2s_0-s)}$, we get 
    $$\|f_t\|_{H^{-s}} \leq \left(C_{\omega}^{s/s_0}+1\right) e^{-s\gamma t/(2s_0-s)} \|f_0\|_{H^s},$$
    completing the proof.
\end{proof}

\revision{
\section{{Compactness of the parabolic solution operator}}\label{app:compact}
In this appendix, we discuss the almost-sure compactness of the solution semigroup $\mathcal T^t_{\kappa,\omega}: f_0 \mapsto f_t$ for the advection-diffusion equation 
\begin{equation}\label{spde:advDif2}
    \dd f_t = \kappa \Delta f_t\, \dd t- \sum_{j = 1}^m \sigma_j \cdot \nabla f_t \circ \dd W_t^j
\end{equation}
in the parabolic case $\kappa>0$. 

\begin{proposition}
   Let $\kappa >0$. For almost every $\omega \in \Omega$, we have that for every $t>0$ the operator $\mathcal T_{\kappa,\omega}^t:L^2(\mathbb T^d) \to L^2(\mathbb T^d)$ is compact.
\end{proposition}
This proposition is used to apply the Multiplicative Ergodic Theorem for compact operators \cite{ruelle1982characteristic} in Section \ref{sec:white}. 
\begin{proof}
    Recall that $\varphi_\omega^{t_0,t}$ denotes the stochastic flow induced by the SDE \eqref{sde:Lagrangian} as introduced in Appendix \ref{app:mixing rates}.
    For a solution $f_t$ of \eqref{spde:advDif2}, let $\tilde f_t$ be given by
    $$\tilde f_t = f_t \circ \varphi^{0,t}_\omega.$$
    An elementary computation shows that $\tilde f_t$ is the solution to the random heat equation
    $$\partial_t \tilde f_t = \kappa \nabla \cdot (a_{\omega,t} \nabla \tilde f_t),$$
     with $\tilde f_0 = f_0$, where $a_{\omega,t}: \mathbb T^d \to \mathbb R^{d\times d}$ is given by
    $$a_{\omega,t}(\mathbf x) \eqqcolon (\mathrm D_{\mathbf x} \varphi_\omega^{0,t}(\mathbf x))^{-1}(\mathrm D_{\mathbf x} \varphi_\omega^{0,t}(\mathbf x))^{-T}.$$
    Here $(\mathrm D_{\mathbf x} \varphi_\omega^{0,t}(\mathbf x))^{-T}$ denotes the inverse transpose of the Jacobian. 
    Since $\varphi$ is a continuous flow of diffeomorphisms, $a_{\omega,t}$ is uniformly elliptic almost surely. Hence, the solution operator $\tilde {\mathcal T}^t_{\kappa,\omega}: L^2(\mathbb T^d)\to L^2(\mathbb T^d), \tilde f_0 \mapsto \tilde f_t$ is compact almost surely. Furthermore, since $f_t = \tilde f_t \circ (\varphi_\omega^{0,t})^{-1}$, the solution operator $\mathcal T_{\kappa,\omega}^t$ can be expressed as the composition $\mathcal T_{\kappa,\omega}^t=\mathcal T_{0,\omega}^t\circ \tilde{\mathcal T}_{\kappa,\omega}^t$, see \eqref{eq:ftConstruction}. Since $\mathcal T_{0,\omega}^t$ is bounded, $\mathcal T_{\kappa,\omega}^t$ is compact as well.
    \end{proof}
}

\section*{Acknowledgements}
The authors thank \revision{Tommaso} Rosati for insightful discussions on this topic. RC thanks Alex Blumenthal and Sam Punshon-Smith for their hospitality during his visits, which initiated the work on this project. \revision{The authors thank the anonymous reviewers for their insightful suggestions which greatly improved the structure of the proofs.}
The work of DC is funded by the Deutsche Forschungsgemeinschaft (DFG, German Research Foundation)
- SPP 2410 Hyperbolic Balance Laws in Fluid Mechanics: Complexity, Scales, Randomness (CoScaRa).
RC acknowledges the funding by the DFG-Collaborative Research Center 1114 Scaling Cascades in Complex Systems, project no. 235221301, A08 Characterization and prediction of quasi-stationary atmospheric states.

\section*{Statements and Declarations}
The authors declare no financial or non-financial competing interests directly or indirectly related to this work. There is no data associated with the results.

\bibliographystyle{alpha}
\bibliography{bibliography.bib}

\end{document}